\documentclass{rmmcart}

     \usepackage{amsmath}
\usepackage{amssymb}
\usepackage{amsthm}
\usepackage{enumerate}
\usepackage[mathscr]{eucal}
\usepackage{eqlist}

\theoremstyle{plain}
\newtheorem{theorem}{Theorem}[section]

\newtheorem{claim}{Claim}
\newtheorem{cla}{Claim}

\theoremstyle{definition}

\numberwithin{equation}{section}

     \title[Non-linear $\ast$-Jordan triple derivation]{Non-linear $\ast$-Jordan triple derivation on prime $\ast$-algebras}

      \author{V. Darvish}
     \address{School of Mathematics and Statistics, Nanjing University of Information Science and Technology,
Nanjing, China.}
     \email{vahid.darvish@mail.com}

    \author{M. Nouri}\author{M. Razeghi}\author{A. Taghavi}
     \address{Department of Mathematics, Faculty of Mathematical
Sciences, University of Mazandaran,
Babolsar, Iran.}
     \email{mojtaba.nori2010@gmail.com}
     
     \email{razeghi.mehran19@yahoo.com}
         \email{taghavi@umz.ac.ir}

     \date{Month, Day, Year}


\begin{document}
     
     \begin{abstract}
     Let $\mathcal{A}$ be a prime $\ast$-algebra and 
 $\Phi$ preserves triple
$\ast$-Jordan derivation on $\mathcal{A}$, that is, for every $A,B \in
\mathcal{A}$,
$$\Phi(A\diamond  B \diamond C)=\Phi(A)\diamond B\diamond  C+A\diamond \Phi(B)\diamond C+A\diamond  B\diamond \Phi(C)$$  
  where $A\diamond  B = AB +  BA^{\ast}$ then $\Phi$ is additive. Moreover, if $\Phi(\alpha I)$ is self-adjoint for $\alpha\in\{1,i\}$ then $\Phi$ is a $\ast$-derivation.
     \end{abstract}
     \maketitle
\section{Introduction}\label{intro}

Let $\mathcal{R}$ be a $*$-ring. For $A,B\in\mathcal{R}$, denoted by
$A\diamond B=AB+BA^{*}$ and $[A,B]_{*}=AB-BA^{*}$, which are $\ast$-Jordan product and $\ast$-Lie product, respectively. These products are found playing a
more and more important role in some research topics, and its study
has recently attracted many author's attention (for example, see
\cite{cui,li2,mol,taghavi}).

Let define $\lambda$-Jordan $\ast$-product by $A\diamond_{\lambda}B =
AB + \lambda BA^{\ast}$. We say the map $\Phi$ with property of
$\Phi(A\diamond_{\lambda} B) = \Phi(A)\diamond_{\lambda}B +  A
\diamond_{\lambda}\Phi(B)$ is a $\lambda$-Jordan $\ast$-derivation
map. It is clear that for $\lambda = - 1$ and $\lambda = 1,$ the $\lambda
$-Jordan $\ast$-derivation map is a $\ast$-Lie derivation and
$\ast$-Jordan derivation, respectively \cite{bai}. We should mention here whenever  we say $\Phi$ preserves derivation, it means $\Phi(AB)=\Phi(A)B+A\Phi(B)$.\\ 
Recently, Yu and Zhang in \cite{yu}
 proved that every non-linear $\ast$-Lie derivation  from a factor von Neumann algebra into itself is an additive $\ast$-derivation. Also,
 Li, Lu and Fang
in \cite{li} have investigated a non-linear $\lambda$-Jordan
$\ast$-derivation. They showed that if
$\mathcal{A}\subseteq\mathcal{B(H)}$ is a von Neumann algebra
without central abelian projections and $\lambda$ is a non-zero scaler,
then $\Phi:\mathcal{A} \longrightarrow \mathcal{B(H)}$ is a
non-linear $\lambda$-Jordan $\ast$-derivation if and only if $\Phi$ is
an additive $\ast$-derivation.

 In \cite{taghavi2} we showed that $\ast$-Jordan derivation map (i.e., $\phi(A\diamond B) = \phi(A) \diamond  B + A\diamond  \phi( B)$)  on every
factor von Neumann algebra $\mathcal{A}
 \subseteq \mathcal{B(H)}$  is
additive $\ast$-derivation. Moreover, in \cite{tvsl} we reduced the assumptions on the map $\phi$.

The authors of \cite{li23} introduced the concept of Lie triple derivations. A map $\Phi:\mathcal{A}\to\mathcal{A}$ is a nonlinear skew Lie triple derivations if $\Phi([[A,B]_{\ast},C]_{\ast})=[[\Phi(A),B]_{\ast},C]_{\ast}+[[A,\Phi(B)]_{\ast},C]_{\ast}+[[A,B]_{\ast},\Phi(C)]_{\ast}$ for all $A, B, C\in\mathcal{A}$ where $[A,B]_{\ast}=AB-BA^{*}$. They showed that if $\Phi$ preserves the above characterizations on factor von Neumann algebras then $\Phi$ is additive $\ast$-derivation.

Also, in \cite{manv} we considered a map $\Phi$ on prime $\ast$-algebra $\mathcal{A}$ which holds in the following conditions
$$\Phi(A\diamond_{\lambda} B \diamond_{\lambda}C)=\Phi(A)\diamond_{\lambda}B\diamond_{\lambda} C+A\diamond_{\lambda}\Phi(B)\diamond_{\lambda}C)+A\diamond_{\lambda} B\diamond_{\lambda}\Phi(C)$$  
  where $A\diamond_{\lambda} B = AB + \lambda BA^{\ast}$ such that a complex scalar $|\lambda|\neq 0, 1$, then $\Phi$ is additive. Also, if $\Phi(I)$ is self-adjoint then $\Phi$ is $\ast$-derivation.

In this paper inspired by the results above, we consider a map $\Phi$ on prime $\ast$-algebra $\mathcal{A}$ which holds in the following conditions
$$\Phi(A\diamond  B \diamond C)=\Phi(A)\diamond B\diamond  C+A\diamond \Phi(B)\diamond C)+A\diamond  B\diamond \Phi(C)$$  
  where $A\diamond  B = AB + BA^{\ast}$  then $\Phi$ is additive. Also, if $\Phi(\alpha I)$ is self-adjoint for $\alpha\in\{1,i\}$ then $\Phi$ is $\ast$-derivation.

We say that $\mathcal{A}$ is
prime, that is, for $A,B \in \mathcal{A}$ if $A\mathcal{A}B =
\lbrace0\rbrace,$ then $A = 0$ or $B = 0$.

\section{Main Results}

Our first theorem is as follows:
\begin{theorem}
Let $\mathcal{A}$ be a prime $\ast$-algebra. Then the map $\Phi:\mathcal{A}\to \mathcal{A}$ satisfies in the following condition 
\begin{equation}\label{sh1}
\Phi(A\diamond  B \diamond C)=\Phi(A)\diamond B\diamond  C+A\diamond \Phi(B)\diamond C+A\diamond  B\diamond \Phi(C)
\end{equation}
for all $A,B\in\mathcal{A}$ where $A\diamond B=AB+BA^{*}$, is additive.
\end{theorem}
\textbf{Proof.}
Let $P_{1}$ be a nontrivial
projection in $\mathcal{A}$ and $P_{2}=I_{\mathcal{A}}-P_{1}$.
Denote $\mathcal{A}_{ij}=P_{i}\mathcal{A}P_{j},\ i,j=1,2,$ then
$\mathcal{A}=\sum_{i,j=1}^{2}\mathcal{A}_{ij}$. For every
$A\in\mathcal{A}$ we may write $A=A_{11}+A_{12}+A_{21}+A_{22}$. In
all that follow, when we write $A_{ij}$, it indicates that
$A_{ij}\in\mathcal{A}_{ij}$. For showing additivity of $\Phi$ on $\mathcal{A}$, we  use above
partition of $\mathcal{A}$ and give some claims that prove $\Phi$ is
additive on each $\mathcal{A}_{ij}, \ i,j=1,2$.

We prove the above theorem by several claims.
\begin{claim}
We show that $\Phi(0)=0$.
\end{claim}
This claim is easy to prove.

\begin{claim}\label{F4}
For each $A_{12} \in \mathcal{A}_{12}$ and $A_{21} \in \mathcal{A}_{21}$ we have 
$$\Phi(A_{12}+A_{21})=\Phi(A_{12})+\Phi(A_{21}).$$
\end{claim}
We show that 
$$T=\Phi(A_{12}+A_{21})-\Phi(A_{12})-\Phi(A_{21})=0.$$
We can write that
\begin{eqnarray*}
&&\Phi\left(\frac{I}{2}\right)\diamond (P_{1}-P_{2})\diamond (A_{12}+A_{21})+\frac{I}{2}\diamond \Phi(P_{1}-P_{2})\diamond (A_{12}+A_{21})\\
&&+\frac{I}{2}\diamond (P_{1}-P_{2})\diamond \Phi(A_{12}+A_{21})\\
&&=\Phi\left(\frac{I}{2}\diamond (P_{1}-P_{2})\diamond (A_{12}+A_{21})\right)\\
&&=\Phi\left(\frac{I}{2}\diamond (P_{1}-P_{2})\diamond A_{12}\right)+\Phi\left(\frac{I}{2}\diamond (P_{1}-P_{2})\diamond A_{21}\right)\\
&&= \Phi\left(\frac{I}{2}\right)\diamond (P_{1}-P_{2})\diamond (A_{12}+A_{21})+\frac{I}{2}\diamond \Phi(P_{1}-P_{2})\diamond (A_{12}+A_{21})\\
&&+\frac{I}{2}\diamond (P_{1}-P_{2})\diamond (\Phi(A_{12})+\Phi(A_{21})).
\end{eqnarray*}
So, we have
$$\frac{I}{2}\diamond (P_{1}-P_{2})\diamond T=0.$$
Since $T=T_{11}+T_{12}+T_{21}+T_{22}$ then 
$T_{11}-T_{22}=0$, therefore $T_{11}=T_{22}$.\\
On the other hand, we obtain
\begin{eqnarray*}
&&\Phi(P_{1})\diamond (A_{12}+A_{21})\diamond P_{1}+P_{1}\diamond \Phi(A_{12}+A_{21})\diamond P_{1}\\
&&+P_{1}\diamond (A_{12}+A_{21})\diamond \Phi(P_{1})\\
&&=\Phi(P_{1}\diamond A_{12}+A_{21}\diamond P_{1})\\
&&=\Phi(P_{1}\diamond A_{12}\diamond P_{1})+\Phi(P_{1}\diamond A_{21}\diamond P_{1})\\
&&=\Phi(P_{1})\diamond (A_{12}+A_{21})\diamond P_{1}\\
&&+P_{1}\diamond (\Phi(A_{12})+\Phi(A_{21}))\diamond P_{1})+P_{1}\diamond (A_{12}+A_{21})\diamond \Phi(P_{1}).
\end{eqnarray*}
So, we have
$$P_{1}\diamond T\diamond P_{1}=0.$$
Since $T=T_{11}+T_{12}+T_{21}+T_{22}$ then $$T_{21}+T_{21}^{*}=0$$
it follows that $T_{21}=0$.\\
In a similar way for $P_{2}$ we can obtain
$$P_{2}\diamond T\diamond P_{2}=0.$$
So, $T_{12}=0$.

\begin{claim}\label{F5}
For each $A_{11} \in \mathcal{A}_{11}$, $A_{12} \in \mathcal{A}_{12}$, $A_{21} \in \mathcal{A}_{21}$, we have 
$$\Phi(A_{11}+A_{12}+A_{21})=\Phi(A_{11})+\Phi(A_{12})+\Phi(A_{21}).$$
\end{claim}
We show that for  $T$ in $\mathcal{A}$ the following holds
\begin{equation}\label{tad4}
T=\Phi(A_{11}+A_{12}+A_{21})-\Phi(A_{11})-\Phi(A_{12})-\Phi(A_{21})=0.
\end{equation}
We can write 
\begin{eqnarray*}
&&\Phi\left(\frac{I}{2}\right)\diamond (P_{1}-P_{2})\diamond (A_{11}+A_{12}+A_{21})+\frac{I}{2}\diamond \Phi(P_{1}-P_{2})\diamond (A_{11}+A_{12}+A_{21})\\
&&+\frac{I}{2}\diamond (P_{1}-P_{2})\diamond \Phi(A_{11}+A_{12}+A_{21})=\Phi\left(\frac{I}{2}\diamond (P_{1}-P_{2})\diamond (A_{11}+A_{12}+A_{21})\right)\\
&&=\Phi\left(\frac{I}{2}\diamond (P_{1}-P_{2})\diamond A_{12}\right)+\Phi\left(\frac{I}{2}\diamond (P_{1}-P_{2})\diamond A_{21}\right)\\
&&+\Phi\left(\frac{I}{2}\diamond (P_{1}-P_{2})\diamond A_{11}\right)=\Phi\left(\frac{I}{2}\right)\diamond (P_{1}-P_{2})\diamond (A_{11}+A_{12}+A_{21})\\
&&+\frac{I}{2}\diamond \Phi(P_{1}-P_{2})\diamond (A_{11}+A_{12}+A_{21})\\
&&+\frac{I}{2}\diamond (P_{1}-P_{2})\diamond (\Phi(A_{11})+\phi(A_{12})+\Phi(A_{21})).
\end{eqnarray*}
Then, we have 
$$\frac{I}{2}\diamond (P_{1}-P_{2})\diamond T=0.$$
Since $T=T_{11}+T_{12}+T_{21}+T_{22}$ we obtain $T_{11}-T_{22}=0$ or $T_{11}=T_{22}=0$.\\
By Claim \ref{F4} we have
\begin{eqnarray*}
&&\Phi(P_{2})\diamond \frac{I}{2}\diamond (A_{11}+A_{12}+A_{21})+P_{2}\diamond \Phi\left(\frac{I}{2}\right)\diamond (A_{11}+A_{12}+A_{21})\\
&&+P_{2}\diamond \frac{I}{2}\diamond \Phi(A_{11}+A_{12}+A_{21})\\
&&=\Phi\left(P_{2}\diamond \frac{I}{2}\diamond (A_{11}+A_{12}+A_{21})\right)\\
&&=\Phi(P_{2}\diamond \frac{I}{2}\diamond A_{11})+\Phi(P_{2}\diamond \frac{I}{2} \diamond (A_{12}+A_{21}))\\
&&=\Phi(P_{2}\diamond \frac{I}{2} \diamond A_{11})+\Phi(P_{2}\diamond \frac{I}{2} \diamond A_{12})+\Phi(P_{2}\diamond \frac{I}{2} \diamond A_{21})\\
&&=\Phi(P_{2})\diamond \frac{I}{2}\diamond (A_{11}+A_{12}+A_{21})+P_{2}\diamond \Phi(\frac{I}{2})\diamond (A_{11}+A_{12}+A_{21})\\
&&+P_{2}\diamond \frac{I}{2} \diamond (\Phi(A_{11})+\Phi(A_{12})+\Phi(A_{21}))
\end{eqnarray*}
So, $$P_{2}\diamond \frac{I}{2} \diamond T=0.$$
Since $T=T_{11}+T_{12}+T_{21}+T_{22}$, we obtain
$$T_{12}+T_{21}+2T_{22}=0.$$
Therefore $T_{12}=T_{21}=0$.

\begin{claim}\label{F6}
For each $A_{11} \in \mathcal{A}_{11}$, $A_{12} \in \mathcal{A}_{12}$, $A_{21} \in \mathcal{A}_{21}$, $A_{22}\in\mathcal{A}_{22}$ we have 
$$\Phi(A_{11}+A_{12}+A_{21}+A_{22})=\Phi(A_{11})+\Phi(A_{12})+\Phi(A_{21})+\Phi(A_{22}).$$
\end{claim}
We show that for  $T$ in $\mathcal{A}$ the following holds
\begin{equation}\label{tad4}
T=\Phi(A_{11}+A_{12}+A_{21}+A_{22})-\Phi(A_{11})-\Phi(A_{12})-\Phi(A_{21}) -\Phi(A_{22})=0.
\end{equation}
From Claim \ref{F5} We can write 
\begin{eqnarray*}
&&\Phi(P_{1})\diamond \frac{I}{2}\diamond (A_{11}+A_{12}+A_{21}+A_{22})+P_{1}\diamond \Phi(\frac{I}{2})\diamond (A_{11}+A_{12}+A_{21}+A_{22})\\
&&+P_{1}\diamond \frac{I}{2} \diamond \Phi(A_{11}+A_{12}+A_{21}+A_{22})\\
&&=\Phi(P_{1}\diamond \frac{I}{2} \diamond (A_{11}+A_{12}+A_{21}+A_{22}))\\
&&=\Phi(P_{1}\diamond \frac{I}{2} \diamond (A_{11}+A_{12}+A_{21}))+\Phi(P_{1}\diamond A_{22})\\
&&=\Phi(P_{1}\diamond \frac{I}{2} \diamond A_{11})+\Phi(P_{1}\diamond \frac{I}{2} \diamond A_{12})+\Phi(P_{1}\diamond \frac{I}{2} \diamond A_{21})+\Phi(P_{1}\diamond \frac{I}{2} \diamond A_{22})\\
&&=\Phi(P_{1})\diamond \frac{I}{2} \diamond (A_{11}+A_{12}+A_{21}+A_{22})+P_{1}\diamond \Phi(\frac{I}{2})\diamond (A_{11}+A_{12}+A_{21}+A_{22})\\
&&+P_{1}\diamond \frac{I}{2}\diamond (\Phi(A_{11})+\Phi(A_{12})+\Phi(A_{21})+\Phi(A_{22})).
\end{eqnarray*}
So,
$$P_{1}\diamond \frac{I}{2}\diamond T=0.$$
Since
$T=T_{11}+T_{12}+T_{21}+T_{22}$ then $T_{11}=T_{21}=T_{22}=0$. Similarly, we can show that $T_{22}=0$.

\begin{claim}\label{F7}
For each $A_{ij},B_{ij} \in \mathcal{A}_{i}$ such that $i\neq j$, we have 
$$\Phi(A_{ij}+B_{ij})=\Phi(A_{ij})+\Phi(B_{ij}).$$
\end{claim}
It is easy to check that 
$$\frac{I}{2}\diamond (P_{i}+A_{ij})\diamond (P_{j}+B_{ij})=A_{ij}+B_{ij}+A_{ij}^{*}+B_{ij}A_{ij}^{*}.$$
From Claim \ref{F6} we have
\begin{eqnarray*}
&&\Phi(A_{ij}+B_{ij})+\Phi(A_{ij}^{*})+\Phi(B_{ij}A_{ij}^{*})=\Phi\left(\frac{I}{2}\diamond (P_{i}+A_{ij})\diamond (P_{j}+B_{ij})\right)\\
&&=\Phi\left(\frac{I}{2}\right)\diamond (P_{i}+A_{ij})\diamond (P_{j}+B_{ij})+\frac{I}{2}\diamond \Phi(P_{i}+A_{ij})\diamond (P_{j}+B_{ij})\\
&&+\frac{I}{2}\diamond (P_{i}+A_{ij})\diamond \Phi(P_{j}+B_{ij})\\
&&=\Phi\left(\frac{I}{2}\right)\diamond (P_{i}+A_{ij})\diamond (P_{j}+B_{ij})+\frac{I}{2}\diamond (\Phi(P_{i})+\Phi(A_{ij}))\diamond (P_{j}+B_{ij})\\
&&+\frac{I}{2}\diamond (P_{i}+A_{ij})\diamond (\Phi(P_{j})+\Phi(B_{ij}))\\
&&=\Phi\left(\frac{I}{2}\diamond A_{ij}\diamond B_{ij}\right)+\Phi\left(\frac{I}{2}\diamond P_{i}\diamond P_{j}\right)+\Phi\left(\frac{I}{2}\diamond P_{i}\diamond B_{ij}\right)\\
&&+\Phi\left(\frac{I}{2}\diamond A_{ij}\diamond P_{j}\right)\\
&&=\Phi(B_{ij})+\Phi(A_{ij}+A_{ij}^{*})+\Phi(B_{ij}A_{ij}^{*}).
\end{eqnarray*}
So, 
$$\Phi(A_{ij}+B_{ij})=\Phi(A_{ij})+\Phi(B_{ij}).$$
\begin{claim}\label{F8}
For each $A_{ii},B_{ii} \in \mathcal{A}_{ii}$ such that $1\leq i \leq 2$, we have 
$$\Phi(A_{ii}+B_{ii})=\Phi(A_{ii})+\Phi(B_{ii}).$$
\end{claim}
We show that 
$$T=\Phi(A_{ii}+B_{ii})-\Phi(A_{ii})-\Phi(B_{ii})=0.$$
We can write
\begin{eqnarray*}
&&\Phi(P_{j})\diamond \frac{I}{2}\diamond (A_{ii}+B_{ii})+P_{j}\diamond \Phi\left(\frac{I}{2}\right)\diamond (A_{ii}+B_{ii})+P_{j}\diamond \frac{I}{2}\diamond \Phi(A_{ii}+B_{ii})\\
&&=\Phi\left(P_{j}\diamond \frac{I}{2}\diamond (A_{ii}+B_{ii})\right)\\
&&=\Phi\left(P_{j}\diamond \frac{I}{2}\diamond A_{ii}\right)+\Phi(P_{j}\diamond \frac{I}{2}\diamond B_{ii})\\
&&=\Phi(P_{j})\diamond \frac{I}{2}\diamond (A_{ii}+B_{ii})+P_{j}\diamond \Phi(\frac{I}{2})\diamond (A_{ii}+B_{ii})\\
&&+P_{j}\diamond \frac{I}{2}\diamond (\Phi(A_{ii})+\Phi(B_{ii})).
\end{eqnarray*}
Therefore, 
$$P_{j}\diamond \frac{I}{2}\diamond T=0.$$
Since 
$T=T_{11}+T_{12}+T_{21}+T_{22}$ we have $T_{jj}=T_{ji}=T_{ij}=0$.\\
From Claim \ref{F7} for every $C_{ij}\in\mathcal{A}_{ij}$ we have
\begin{eqnarray*}
&&\Phi(P_{i})\diamond (A_{ii}+B_{ii})\diamond C_{ij}+P_{i}\diamond \Phi(A_{ii}+B_{ii})\diamond C_{ij}+P_{i}\diamond (A_{ii}+B_{ii})\diamond \Phi(C_{ij})\\
&&=\Phi(P_{i}\diamond (A_{ii}+B_{ii})\diamond C_{ij})\\
&&=\Phi(P_{i}\diamond A_{ii}\diamond  C_{ij})+\Phi(P_{i}\diamond B_{ii}\diamond C_{ij})\\
&&=\Phi(P_{i})\diamond (A_{ii}+B_{ii})\diamond C_{ij}+P_{i}\diamond (\Phi(A_{ii})+\Phi(B_{ii}))\diamond C_{ij}\\
&&+P_{i}\diamond (A_{ii}+B_{ii})\diamond \Phi(C_{ij}).
\end{eqnarray*}
So, 
$$P_{i}\diamond T\diamond C_{ij}=0.$$
By primeness and since $T=T_{11}+T_{12}+T_{21}+T_{22}$, by primeness we obtain $T_{ii}=0$.\\
Hence, the additivity of $\Phi$ comes from the above claims.
\qed

In the rest of this paper we prove that $\Phi$ is $\ast$-derivation.
\begin{theorem}
Let $\mathcal{A}$ be a prime $\ast$-algebra. Let $\Phi:A\to A$ satisfies in the following condition 
\begin{equation}\label{sh1}
\Phi(A\diamond  B \diamond C)=\Phi(A)\diamond B\diamond  C+A\diamond \Phi(B)\diamond C)+A\diamond  B\diamond \Phi(C)
\end{equation}
for all $A,B\in\mathcal{A}$ where $A\diamond B=AB+BA^{*}$. If $\Phi(\alpha I)$ is self-adjoint for $\alpha\in\{1,i\}$ then $\Phi$ is $\ast$-derivation.
\end{theorem}
\textbf{Proof.}
We present the proof by several claims.
\begin{cla}\label{ai}
If $\Phi(I)$ is self-adjoint then $\Phi(I)=0$.
\end{cla}
One can easily show that 
$$\Phi(I\diamond I\diamond I)=\Phi(I)\diamond I\diamond I+I\diamond \Phi(I)\diamond I+I\diamond I\diamond \Phi(I).$$
So, $\Phi(I)=0$.
\begin{cla}\label{aii}
If $\Phi(iI)$ is self-adjoint then $\Phi(iI)=0$.
\end{cla}
It is easy to check that 
$$\Phi(iI\diamond I\diamond I)=\Phi(iI)\diamond I\diamond I+iI\diamond \Phi(I)\diamond I+iI\diamond I\diamond \Phi(I)=0.$$
We obtain
$$\Phi(iI)\diamond I\diamond  I=0.$$
So, $\Phi(iI)+\Phi(iI)^{*}=0$. Since $\Phi(iI)$ is self-adjoint then $\Phi(iI)=0$.

\begin{cla}
$\Phi$ preserves star.
\end{cla}
Since $\Phi(I)=0$, then
$$\Phi(I\diamond A\diamond I)=I\diamond \Phi(A)\diamond I.$$
Therefore
$$2\Phi(A+A^{*})=2\Phi(A)+2\Phi(A)^{*}.$$
So, we obtain
$$\Phi(A^{*})=\Phi(A)^{*}.$$
\begin{cla}
We show that $\Phi(iA)=i\Phi(A)$ for every $A\in\mathcal{A}$.
\end{cla}
Let $T$ be a self-adjoint member of $\mathcal{A}$. Then it is easy to check that
$$\Phi\left(I\diamond T\diamond iI\right)=\Phi\left(I\diamond I\diamond iT\right).$$
Since $\Phi\left(I\right)=\Phi\left(iI\right)=0$, we obtain
$$I\diamond \Phi(T)\diamond iI=I\diamond I\diamond \Phi(iT).$$
So, we have
$$2i(\Phi(T)+\Phi(T)^{*})=2(\Phi(iT)+\Phi(iT)).$$
Since $\Phi$ preserves star and is additive, for every self-adjoint $T$ we have
$$i\Phi(T)=\Phi(iT).$$
On the other hand, we can write every $A\in\mathcal{A}$ as $A=A_{1}+iA_{2}$ where $A_{1}=\frac{A+A^{*}}{2}$ and $A_{2}=\frac{A-A^{*}}{2i}$.
Finally we have
\begin{eqnarray*}
\Phi(iA)&=&\Phi(iA_{1}-A_{2})\\
&=&\Phi(iA_{1})-\Phi(A_{2})\\
&=&\Phi(iA_{1})+i^{2}\Phi(A_{2})\\
&=&i\Phi(A_{1})+i\Phi(iA_{2})\\
&=&i\Phi(A_{1}+iA_{2})=i\Phi(A).
\end{eqnarray*}
So, $\Phi(iA)=i\Phi(A)$
for all $A\in\mathcal{A}$.

\begin{cla}
$\Phi$ is a derivation.
\end{cla}
For every $A, B\in\mathcal{A}$ we have
\begin{eqnarray*}
2\Phi(AB+BA^{*})&=&\Phi(2AB+2BA^{*})\\
&=&\Phi(I\diamond A\diamond B)\\
&=&\Phi(I)\diamond A\diamond B+I\diamond \Phi(A)\diamond B+I\diamond A\diamond  \Phi(B)\\
&=&2\Phi(B)A^{*}+2\Phi(A)B+2B\Phi(A)^{*}+2A\Phi(B).
\end{eqnarray*}
It follows that
\begin{equation}\label{tam}
\Phi(AB+BA^{*})=\Phi(B)^{*}A+\Phi(A)B+B\Phi(A)^{*}+A\Phi(B).
\end{equation}
From (\ref{tam}) we have 
\begin{eqnarray*}
&&\Phi(AB-BA^{*})=\Phi((iA)(-iB)+(-iB)(iA)^{*})\\
&&=\Phi(iA)(-iB)+(iA)\Phi(-iB)+\Phi(-iB)^{*}(iA)+(-iB)\Phi(iA)^{*}\ \ \text{By (\ref{tam})}\\
&&=\Phi(A)B+A\Phi(B)-\Phi(B)A^{*}-B\Phi(A)^{*}.
\end{eqnarray*}
So, 
\begin{equation}\label{tam2}
\Phi(AB-BA^{*})=\Phi(A)B+A\Phi(B)-\Phi(B)A^{*}-B\Phi(A)^{*}.
\end{equation}
From (\ref{tam}) and (\ref{tam2}) we obtain
$$\Phi(AB)=\Phi(A)B+A\Phi(B).$$
\qed
\paragraph*{}
\paragraph*{}

\textbf{Acknowledgments:} 
The first author is supported by the Talented Young Scientist Program of Ministry of Science and Technology
of China (Iran-19-001).

     \end{document}